\magnification=\magstephalf
\input amstex
\documentstyle{amsppt}
\pagewidth{34pc} \pageheight{47pc}
\hcorrection{.125truein}\vcorrection{.125truein}
\loadbold
\loadeusm 
 
\TagsAsMath
\def\memo#1{}
   \def\Aut{\operatorname{Aut}} 
   \def\card{\operatorname{card}}
   \def\Mir{\operatorname{Mir}}   
   \def\a{\alpha} \undefine\b \def\b{\beta} \def\d{\delta}
   \def\e{\varepsilon} \def\g{\gamma} 
   \def\th{\theta} 
   \def\phi{\varphi} \def\r{\varrho} \def\s{\sigma}  
   \def\D{\Delta}
   \def\G{\Gamma}

     \def\after{\circ}
     \def\Bd{\partial} 
     \def\char#1{\bigl\langle#1\bigr\rangle}
     \def\Nb#1#2{N_#1(#2)}
     \def\C{\Bbb C}  
     \def\R{\Bbb R} \def\Z{\Bbb Z}
     \def\h#1#2{h^{\sssize (#1)}_{#2}} 
     \def\sub{\subset} \redefine\emptyset{\varnothing}
     \font\bxtwelve=cmbx12            
          \def\star{\botsmash{\lower.75ex\hbox{\text{\bxtwelve*}}}}
     \def\plumb#1{\ast_#1} 
     \def\PLUMB{\mathop{\star}\nolimits}
\hyphenation{arbor-es-cent 
homo-gen-e-ous homo-topy
Worces-ter}
\topmatter
\thanks
Research partially supported by CNRS 
and by NSF grant DMS-9504832.
\endthanks
\title 
Quasipositive pretzels
\endtitle
\rightheadtext{Quasipositive pretzels}
\leftheadtext{Lee Rudolph}
\author 
Lee Rudolph 
\endauthor
\address 
Department of Mathematics, Clark University,
Worcester MA 01610, USA
\endaddress
\email 
lrudolph\@black.clarku.edu
\endemail
\abstract
A necessary and sufficient condition 
for an oriented pretzel surface to be quasipositive
yields an estimate for the slice genus of the boundary
of an arbitrary oriented pretzel surface.
\endabstract
\endtopmatter
\document

\head \S1. Introduction\endhead

A {\it braidzel} (see \S2) is a generalized pretzel surface,
with braiding data supplementing the twisting data 
that specifies an ordinary pretzel surface (a braidzel with 
trivial braiding).
{\it Quasipositive\/} Seifert surfaces (see \S3) 
can be characterized as those generated from $D^2\sub S^3$
by plumbing positive Hopf annuli $A(O,-1)$ and 
passing to $\pi_1$-injective subsurfaces; 
they are analogues in $S^3$ (and, in a sense, 
special cases) of pieces of complex plane curve in 
$D^4\sub\C^2$.

\proclaim{Theorem} The oriented pretzel surface 
$P(t_1,\dots,t_k)$ is quasipositive if and only if
the even integer $t_i+t_j$ is less than $0$ for 
$1\le i<j\le k$.
\endproclaim

This result (announced in \cite{7} for $k=3$ and $t_i$ odd)
is deduced in \S4 from results on oriented braidzels.
An estimate for the slice genus of the boundary of an
oriented pretzel surface is given in \S5.  
Other applications will appear elsewhere \cite{0}.

I thank Michel Boileau for his help and encouragement with this
and related matters over the years, and Universit\'e Paul Sabatier 
and CNRS for their support.

\head \S2. Braidzels\endhead

Let 
$E_k(F):=\{\{w_1,\dots,w_k\}\sub F:0\not=\prod\nolimits_{i\ne j}w_i-w_j\}$
be the configuration space of unordered $k$-tuples of elements of the
field $F$.
A {\it $k$-string braid group\/} is the fundamental group
$B_X:=\pi_1(E_k(\C);X)$ with arbitrary basepoint $X\in E_k(\C)$;
of course, every $k$-string braid group is isomorphic to the 
standard $k$-string braid group $B_k:=B_{\{1,\dots,k\}}$
but it is very convenient to allow more general basepoints.  
Typically (for $k>2$) no isomorphism
between distinct $k$-string braid groups has much claim 
to be called canonical, but (because $E_k(\R)$ is contractible,
as is easily seen) if $X_0, X_1\in E_k(\R)$ then 
there is a canonical homotopy class of paths in $E_k(\C)$ 
from $X_0$ to $X_1$ (namely, those that stay in $E_k(\R)$), 
and thus a canonical isomorphism from $B_{X_0}$ to $B_{X_1}$.
Consequently, any path $p:([0,1];0,1)\to (E_k(\C);X_0,X_1)$ 
represents a well-defined braid $\b(p)\in B_k$---in particular, 
for any $X\in E_k(\R)$ there are well-defined braids 
$\D_X$, $\d_X$ and $\r_X$ in $B_X$ which are in canonical 
correspondance with 
$\D_k:=(\s_{k-1})(\s_{k-2}\s_{k-1})\dotsm(\s_1\s_2\dotsm\s_{k-1})$,
$\d_k:=\s_{k-1}\s_{k-2}\dotsm\s_1$, and 
$\r_k:=\D_k\d_k\D_k^{-1}=\s_1\s_2\dotsm\s_{k-1}$ 
in $B_k$.  Let $\g(p):=\{(z,\th): z\in p(\th), \th\in [0,1]\}\supset %
X_0\times\{0\}\cup X_1\times\{1\}$, 
$\G(p):=[\min X_0,\max X_0]\times\{0\}%
\thinspace\cup\thinspace\g(p)\thinspace\cup\thinspace%
[\min X_1,\max X_1]\times\{1\}$ (Figure~1).
For $X'_0\sub X_0$ with $k':=\card X'_0>0$, 
the union of the components 
of $\g(p)$ which intersect $X'_0\times\{0\}$ 
is $\g(p')$ for a suitable 
$p':([0,1];0,1)\to (E_{k'}(\C);X'_0,X'_1)$, and 
$\G(p')\sub\G(p)$.
Note that $\b=\b(p)$ determines $\b(p')=:\b|_{X'_0}$;
in particular, if $X'_0=\{x_i,x_j\}$, $x_i\ne x_j$, 
then $c(x_i,x_j;\b)\in\Z$ is well-defined by 
$\b|_{X'_0}=\s_1^{c(x_i,x_j;\b)}$.
\topinsert
\hskip.25in
\special{em:graph prtzfig1.pcx}
\vspace{.85in}
\botcaption{Figure 1.  
{\rm $\G(p)$, where 
$\b(p)=\s_2\s_1\s_2\s_3\s_4\s_2\s_5\s_3\s_1\in B_6$.}
}
\endcaption
\endinsert

A {\it braidzel\/} is any 
(unoriented, and not necessarily orientable)
$2$-manifold-with-boundary $S\sub\C\times\R\cong\R^3\sub S^3$ 
which has a $(0,1)$-handle decomposition
$$
S=\bigcup_{s\in\{0,1\}}\h0s\cup \bigcup_{x\in X_0}\h1x \tag{2.1}
$$
related as follows to a (piecewise-smooth) path 
$p:([0,1];0,1)\to (E_k;X_0,X_1)$ with $X_0, X_1\in E_k(\R)$:
\roster\widestnumber\item{(2.1.2)}
\item"(2.1.0)" 
$\h0s\sub\R\times\{\th: s+(-1)^s\th\le 0\}$
and $\h0s\cap(\C\times\{s\})\supset [\min X_s,\max X_s]\times\{s\}$;
\item"(2.1.1)" 
the component of $\g(p)$ 
containing $(x,0)$ is a core arc of $\h1x$.
\endroster
The {\it braiding\/} of $S$ is $\b_S:=\b(p)\in B_k$; 
the {\it twisting\/} of $S$ is the function $t_S:X_0\to\Z$ 
such that
\roster\widestnumber\item{(2.1.2)}
\item"(2.1.2)" $\h1x$ makes $t_S(x)$ counterclockwise half-twists 
between $\h1x\cap\h00$ and $\h1x\cap\h01$ (its two attaching arcs). 
\endroster
Clearly $S$ is orientable
if and only if $t_S(x)+t_S(x')\in 2\Z$ for all $x, x'\in X_0$.

Though the path $p$ to which a given braidzel $S$ is related 
as in \thetag{2.1} is never unique, evidently both $\b_S$ 
and $t_S$ (up to an increasing bijection of its domain) 
do depend on $S$ alone, and not on $p$.  
Also evidently, 
every pair $(\b,t)$ is $(\b_S,t_S)$ for a braidzel 
well-defined up to isotopy through braidzels,
and $P(\b|_{X'_0},t|X'_0)$ can be taken to be a subbraidzel
of $P(\b,t)$ for $\emptyset\ne X'_0\sub X_0$.
It is customary for a braidzel $P(o,t)$ with trivial braiding
$o=o^{(k)}:=1_{B_k}$ (and, implicitly, with $X_0=X_1$) to be called 
a {\it pretzel surface}; in one tradition, 
followed below,
the notation for $P(o,t)$ is $P(t_1,\dots,t_k)$,
where $X_0=:\{x_1,\dots,x_k\}$, $x_1<\dots< x_k$,
and $t_i:=t(x_i)$.

Since $\G(p)$ is a spine of $P(\b,t)$ for all $t$,
a diagram of $\G(p)$ with the components of $\g(p)$
labelled according to $t$ may be used to portray
$P(\b,t)$ schematically (Figure~2).  
\topinsert
\hskip.25in
\special{em:graph prtzfig2.pcx}
\vspace{1.25in}
\botcaption{Figure 2.  
{\rm $P(3,-5,-7)$ (in full-dress; schematic).}
}
\endcaption
\endinsert

For $Z\in E_k(\C)$, write 
$\pi:B_Z\to\Aut(Z)$ for the natural homomorphism;
for $z\in Z$, let $\char{z}:Z\to\{0,1\}\sub\Z$ 
denote the characteristic function of $\{z\}$. 
Let $P(\b,t)$ be a braidzel related as in \thetag{2.1} 
to a path $p:([0,1];0,1)\to (E_k;X,Y)$.

\proclaim{2.2. Lemma} Each of the braidzels \rom{(2.2.1--3)}
\rom{(}resp., \rom{(2.2.4--6))} 
is isotopic to $P(\b,t)$ by an isotopy supported in a 
regular neighborhood of $\h00$ \rom{(}resp., $\h01$\rom{)}.
\par\noindent
\rom{(2.2.1)}~$P(\d_X\b,t\after\pi(\d_X)+2\char{\max X})$
\par\noindent
\rom{(2.2.2)}~$P(\r^{-1}_X\b,t\after \pi(\r^{-1}_X)-2\char{\max X})$
\par\noindent
\rom{(2.2.3)}~$P(\D_X\b,t+1)$
\par\noindent
\rom{(2.2.4)}~$P(\b\d_Y,t+2\char{\pi(\b^{-1})(\max Y)})$
\par\noindent
\rom{(2.2.5)}~$P(\b\r_Y^{-1},t-2\char{\pi(\b^{-1})(\max Y)})$
\par\noindent
\rom{(2.2.6)}~$P(\b\D_Y,t+1)$
\endproclaim
\topinsert
\hskip.25in
\special{em:graph prtzfig3.pcx}
\vspace{1.45in}
\botcaption{Figure 3.  
{\rm $\Nb{{S^3}}{{\h00}}$; 
$1$-handles sliding along $\Bd\h00$;
$P(\d_X\b,t\after \pi(\d_X)+2\char{\max X})$;
$P(\d_X\b,t\after \pi(\r^{-1}_X)-2\char{\max X})$.}
}
\endcaption
\endinsert
\topinsert
\hskip.25in
\special{em:graph prtzfig4.pcx}
\vspace{1.45in}
\botcaption{Figure 4.   
{\rm $P(\b,t)$;
$P(\d_X\b,t\after \pi(\d_X)-2\char{\max X})$;
$P(\r_X\b,t\after \pi(\r_X)+2\char{\max X})$;
$P(\D_X\b,t+1)$ (schematics).}
}
\endcaption
\endinsert
\demo{Proof} For (2.2.1--2), the required isotopy begins with a slide
of all the $1$-handles counterclockwise along $\Bd\h00$; next, for
(2.2.1) (resp., (2.2.2)) the $1$-handle originally containing
$(\min X,0)$ is pulled up in back (resp., in front) of $\h00$ 
and acquires a positive (resp., negative) kink; 
finally, the kink is replaced with a half-twist counterlockwise 
(resp., clockwise), bringing the isotoped surface back into braidzel 
position.  For (2.2.3), the required isotopy is a rotation of $\h00$ 
through one half-twist clockwise.  
(See Figures~3 and 4 for (2.2.1--3).)
Similar isotopies work for (2.2.4--6).
\qed
\enddemo
\remark{Remarks} (1)~Of course 
$P(\d_X^{-1}\b,t\after\pi(\d_X^{-1})-2\char{\max X})$,
\dots, $P(\b\D_Y^{-1},t-1)$ are also isotopic to $P(\b,t)$.  
(2)~It is well known (and obvious) 
that the pretzel surfaces $P(t_1,\dots,t_k)$, $P(t_2,\dots,t_k,t_1)$, 
and $P(t_k,t_{k-1},\dots,t_1)$ are mutually isotopic.
This also follows immediately from (2.2) and the preceding
remark, since $o^{(k)}=\d_k o^{(k)}\d_k^{-1}=\D_k o^{(k)}\D_k^{-1}$.
\endremark 

\head \S3. Plumbing and quasipositivity\endhead
The isotopy type of $P(t_1,t_2)$ depends only on 
$t_1+t_2$.  Any Seifert surface isotopic to 
$P(t,t)$ (endowed with either orientation) is denoted
$A(O,t)$, and called a {\it $t$-twisted unknotted}
annulus in $S^3$ (``unknotted'' because its core 
circle---equivalently, either component of its boundary---is
an unknot $O$, and ``$t$-twisted'' because the Seifert
self-linking of $[O]\in H_1(A(O,t);\Z)$ is $t$).
In particular, $A(O,-1)$ is called a {\it positive 
Hopf annulus\/} ({\it Hopf\/} because the oriented link 
$\Bd A(O,-1)$ consists of two fibers of a Hopf fibration 
$S^3\to S^2$, {\it positive\/} because the linking number 
of these fibers is $+1$).

There is an extensive literature on {\it Murasugi sum\/}
(or {\it Stalling plumbing\/}) and its various special cases
({\it annular plumbing\/}, {\it arborescent plumbing\/}, etc\.);
further details are given in \cite{8}\memo{QPV} and \cite{10}\memo{QPVI},
and references cited therein.  For present purposes it suffices to 
recall {\it positive Hopf plumbing} briefly.

Let $\a\sub S$ be a proper arc on a Seifert surface $S\sub S^3$.
Let $C_\a\sub S$ be a regular neighborhood of $\a$
(naturally endowed with the structure of a $4$-gon whose 
edges are alternately arcs on $\Bd S$ and proper arcs parallel 
to $\a$).  Let $D_\a\sub S^3$ be a $3$-cell on the positive side
of $S$ with $D_\a\cap S=C_\a$. Let $A(O,-1)\sub D_\a$ be a
positive Hopf annulus with $A(O,-1)\cap\Bd D_\a=C_\a$ such
that $\Bd A(O,-1)\sub\Bd C_\a$ consists of the two proper arcs
parallel to $\a$.  Say that $S\cup A(O,-1) =: S\plumb\a A(O,-1)$ 
is constructed by {\it plumbing $A(O,-1)$ to $S$ along $\a$\/},
and that $S$ is constructed from $S\plumb\a A(O,-1)$ by {\it deplumbing\/}
$A(O,-1)$.  Up to isotopy, $S\plumb\a A(O,-1)$ depends only 
on $S$ (up to isotopy) and $\a$ (up to isotopy
on $S$); moreover (Figure 5), $S\plumb\a A(O,-1)$ 
is isotopic to $-(-S\plumb\a A(O,-1))$, 
where $-S$ denotes the surface $S$ with orientation 
reversed.
\topinsert
\hskip.25in
\special{em:graph prtzfig5.pcx}
\vspace{.65in}
\botcaption{Figure 5.
{\rm 
An isotopy from
$S\plumb \a A(O,-1)$ to $-(-S\plumb \a A(O,-1))$.
}
}
\endcaption
\endinsert
A subsurface $T$ of a surface $S$ is {\it $\pi_1$-injective\/}
if, whenever a simple closed curve $C\sub T$ bounds a disk on $S$, 
then $C$ already bounds a disk on $T$.  The notation $T\Subset S$
(or $S\Supset T$) 
will indicate that $T$ is a $\pi_1$-injective subsurface of $S$.  
For example, if $S$ is constructed from $S'$ by deplumbing $A(O,-1)$, 
then $S\Subset S'$.

The notion of {\it quasipositivity\/}, 
introduced in \cite{3} and \cite{4},
has been elaborated and explored in a number of later papers.  
In particular, \cite{5}\memo{III} and \cite{8}\memo{V} 
taken together justify the following definition: 
a Seifert surface $S\sub S^3$ is {\it quasipositive\/} if and
only if, up to isotopy, $S$ is constructed from $D^2$
by a sequence of moves, each of which is either plumbing $A(O,-1)$ 
or passing to a $\pi_1$-injective subsurface.  

\remark{Remarks} (1)~In fact, \cite{5}\memo{III} shows that,
if $S$ is quasipositive, then there exists such a sequence of moves
in which all the positive Hopf plumbings come first (and can even
be done simultaneously, to appropriate arcs on $D^2$, thus creating 
a positive Hopf-plumbed {\it basket\/} in the sense of \cite{8}\memo{QPV}, 
\cite{10}\memo{QPVI}), followed by a single passage to a 
$\pi_1$-injective subsurface (called a {\it full\/} subsurface 
in \cite{5}\memo{III}); the present formulation is more convenient.  
(2)~Also, \cite{8}\memo{V} shows that arbitrary Murasugi sums with
quasipositive plumbands (not just positive Hopf plumbings 
$S\plumb \a A(O,-1)$ with $S$ quasipositive)
are quasipositive; that generality is not needed here.  
\endremark

\head \S4. Characterization of quasipositive pretzel Seifert 
surfaces\endhead
Let $P(\b,t)$ be an orientable braidzel (endowed with either
orientation) related as in \thetag{2.1} to a path 
$p:([0,1];0,1)\to (E_k;X,Y)$.

\proclaim{4.1. Lemma} If $k=2$, so $X=\{x_1,x_2\}$ and 
$\b=\s_1^m=\D_2^m$, then $P(\b,t)$ is quasipositive 
if and only if $t(x_1)+t(x_2)<2m$.
\endproclaim\demo{Proof}
If $m=0$, then $P(\b,t)$ is a $(t(x_1)+t(x_2))/2$-twisted 
unknotted annulus, and the conclusion follows from \cite{6}\memo{IV}.
If $|m|>0$, then $P(\b,t)$ is isotopic (by $|m|-1$ applications 
of (2.2.3)) to $P(t(x_1)-m,t(x_2)-m)$.
\qed\enddemo

\proclaim{4.2. Corollary} If $P(\b,t)$ is quasipositive, then 
$t(x_i)+t(x_j)<2c(x_i,x_j;\b)$ for all $x_i, x_j\in X$, $x_i\ne x_j$.
In particular, if $P(t_1,\dots,t_k)$ is quasipositive, then 
$t_i+t_j<0$ for $1\le i<j\le k$.\qed\endproclaim

Call $t:X\to\Z$ {\it nearly negative\/} if
$t\le 0$ and $\card\{x:t(x)=0\}\le 1$.

\proclaim{4.3. Proposition} If an oriented pretzel surface 
has nearly negative twisting, then it is quasipositive.
\endproclaim
\demo{Proof} If $P(t_1,\dots,t_k)$ is an oriented pretzel surface
with nearly negative twisting, 
then its boundary is a {\it positive link\/},
that is, it has a diagram (namely, the full-dress diagram 
for $P(t_1,\dots,t_k)$ in the style of Figure~2) in which 
every crossing is positive.
By \cite{2} or \cite{9}, the application of 
Seifert's original algorithm \cite{11} to 
this diagram produces a quasipositive Seifert surface $S$,
which by inspection is $P(t_1,\dots,t_k)$.
\qed \enddemo
\remark{Remarks} (1)~Of course the boundary of an oriented
pretzel surface with $t\le 0$ is always a positive link, 
but if $\card\{x:t(x)=0\}>1$ then Seifert's algorithm 
produces a surface which is disconnected, and so 
not $P(t_1,\dots,t_k)$.
(2)~For variety and to illustrate other
techniques, here is an alternative proof using the Twist Insertion 
Lemma \cite{3} (Figure~6 indicates a ``coordinate-free'' statement
and proof of this lemma independent of the machinery of braided
surfaces applied in \cite{3}).  If either all $t_i$ are $-1$, or 
one $t_i$ is $0$ and the rest are $-2$, an easy induction on $k$ 
establishes that $P(t_1,\dots,t_k)$ is isotopic to a quasipositive 
Hopf plumbed basket
$$
(\dots(D^2\plumb{{\a_1}}A(O,-1))\plumb{{\a_2}}A(O,-1)\dots)%
\plumb{{\a_{k-1}}}A(O,-1),
$$
for appropriate proper arcs $\a_1,\dots,\a_{k-1}\sub D^2$
(see Figure~7).
In general, if $t_i\le 0$ for $i=1,\dots,k$, and $t_i=0$ for 
at most one $i$, then $P(t_1,\dots,t_k)$ can (depending on
the parity of the $t_i$) be produced from either $P(-1,\dots,1)$ 
or (at least) one of 
$P(0,-2,\dots,-2)$, $P(0,-2,\dots,-2)$, \dots, $P(-2,-2,\dots,0)$,
by introducing extra clockwise full twists into some of the 
$1$-handles, so the Twist Insertion Lemma shows that 
$P(t_1,\dots,t_k)$ is quasipositive.  \qed
\endremark
\topinsert
\hskip.25in
\special{em:graph prtzfig6.pcx}
\vspace{.65in}
\botcaption{Figure 6.
{\rm The regular neighborhood of an arbitrary proper arc
on a quasipositive Seifert surface $S$; 
$S'$, produced by plumbing two positive Hopf annuli to $S$, 
is also quasipositive; so is $S''\Subset S'$;
but $S''$ is isotopic to $S$ with a clockwise full twist inserted.
}
}
\endcaption
\endinsert
\topinsert
\hskip.25in
\special{em:graph prtzfig7.pcx}
\vspace{1.25in}
\botcaption{Figure 7.
{\rm Plumbing arc arrangements for basket presentations
of $P(-1,\dots,-1)$ and $P(0,-2,\dots,-2)$.}
}
\endcaption
\endinsert
\topinsert
\hskip.25in
\special{em:graph prtzfig8.pcx}
\vspace{1.1in}
\botcaption{Figure 8.   
{\rm A braidzel $P(\b_0\s_i,t)$; a Seifert surface $S$ 
isotopic to $P(\b_0\s_i,t)$; a Seifert surface 
$S'\PLUMB A(O,-1)\Supset S$; a braidzel $P(\b_0,t)$
isotopic to $S'$.}
}
\endcaption
\endinsert

Call a braid $\b\in B_k$ {\it non-negative\/} if it can be written
as a word in the standard generators $\s_1,\dots,\s_{k-1}\in B_k$ 
with no negative exponents.  (A non-negative braid in which each 
$\s_i$ appears non-trivially is conventionally called
{\it positive}.)

\proclaim{4.4. Proposition} If $\b$ is non-negative and $t$ is nearly
negative, then $P(\b,t)$ is quasipositive.
\endproclaim
\demo{Proof} The case of trivial braiding $\b=o^{(k)}$ is 4.3.  The 
proof for general non-negative $\b\in B_k$ is established by induction 
on the length of a word for $\b$ in $\s_1,\dots,\s_{k-1}\in B_k$ without
negative exponents.  The proof of the inductive step (illustrated in
Figure~8) consists in showing, by a sequence of handle slides 
supported in a regular neighborhood of $\h01$, that 
$P(\b_0\s_i,t)$ (resp., $P(\b_0,t)$) 
is isotopic to a Seifert surface $S$ (resp., $S'$)
with $S\Subset S'\PLUMB A(O,-1)$.
\qed\enddemo

\proclaim{4.5. Theorem} An oriented pretzel surface $P(t_1,\dots,t_k)$ 
is quasipositive if and only if $t_i+t_j<0$ for $1\le i<j\le k$.
\endproclaim
\demo{Proof} ``Only if'' was proved in 4.2, and ``if'' with the
extra hypothesis that $t$ be nearly negative in 4.3.  If 
$P(t_1,\dots,t_k)$ is an oriented pretzel surface with 
$t_i+t_j<0$ for $1\le i<j\le k$, and $t$ is not nearly
negative, then there is exactly one $i$ with $t_i>0$, and 
$t_j<0$ for $j\ne i$; without loss of generality (by the second 
remark after 2.2), $t_1>0$.  By (2.2.3) and (2.2.2), 
$P(t_1,t_2,\dots,t_k)$ is isotopic to $P(\D_k\r_k^{-1},t')$, where 
$\D_k\r_k^{-1}=(\s_{k-1})(\s_{k-2}\s_{k-1})\dots(\s_2\s_3\dots\s_{k-1})$ 
is non-negative, $t'(1)=t_1-1\ge 0$,
and $t'(j)=t_{k-1-j}+1<0$ for $2\le j\le k$; a sequence of
$t_1$ such moves gives an isotopy from $P(t_1,t_2,\dots,t_k)$ 
to $P(\b,t'')$, where $\b=(\D_k\r_k^{-1})^{t_1}$ is non-negative
and $t''$ is nearly negative.  By 4.4, $P(t_1,t_2,\dots,t_k)$ 
is quasipositive.
\qed\enddemo

\head \S5. Slice estimates\endhead

Recall that, if $L\sub S^3=\Bd D^4$ is an oriented link,
then $\chi_s(L)$ denotes the greatest Euler characteristic 
$\chi(F)$ such that $F\sub D^4$ is a smooth surface with 
$L=\Bd F$ such that no component of $F$ has empty boundary.
If $K$ is a knot, then its 
{\it slice {\rm (or {\it Murasugi\/})\/} genus\/} $g_s(K)$ (the
least genus $g(F)$ for such an $F$ with $K=\Bd F$) equals
$(1-\chi_s(K))/2$.  

The following results are proved in \cite{7} and \cite{8}.

\proclaim{5.1. Proposition} If $Q$ is a quasipositive Seifert
surface, then $\chi_s(\Bd Q)= \chi(Q)$.  If $S$ is any Seifert surface, 
and $Q\sub S$ is quasipositive, then $\chi_s(\Bd S)\le 2\chi(Q)-\chi(S)$.
\qed
\endproclaim
 
Let $P(\b,t)$ be an oriented braidzel related 
to $p:([0,1];0,1)\to (E_k;X,Y)$.

\proclaim{5.2. Proposition} If $\b$ is non-negative 
and $t$ is nearly negative, then $\chi_s(\Bd P(\b,t))=2-k$. 
More generally, for any $X'\sub X$ such that $\b|X'$ is non-negative and 
$t|X'$ is nearly negative, $\chi_s(\Bd P(\b,t))\le 2-2\card X'+k$. 
\qed
\endproclaim
\demo{Proof} Immediate from 4.4 and 5.1. \qed
\enddemo

Let $P(t_1,\dots,t_k)$ be an oriented pretzel surface.

\proclaim{5.3. Proposition} If $t_i+t_j<0$ for $1\le i<j\le k$,
then $\chi_s(\Bd P(t_1,\dots,t_k))=2-k$.  More generally,
for any $(t_1,\dots,t_k)$, $\chi_s(\Bd P(t_1,\dots,t_k))\le %
2+\card\{i:t_i\ge 0\}-\card\{j:t_j < 0\} - \e$,
where $\e=1$ if $\min\{t_i : t_i\ge 0\}+\max\{t_j : t_j<0\}<0$,
$\e=0$ otherwise.
\endproclaim\demo{Proof} Immediate from 4.5 and 5.1.\qed\enddemo

\remark{Remarks} (1)~Let $\Mir: S^3\to S^3$ be an orientation-reversing
diffeomorphism, so for any oriented $M\sub S^3$, $\Mir M$ is the mirror
image of $M$.  Since $\chi_s(L)=\chi_s(\Mir L)$ for any $L$, 
the application of 5.3 to $P(-t_k,\dots,-t_1)=\Mir P(t_1,\dots,t_k)$
gives another upper bound on $\chi_s(\Bd P(t_1,\dots,t_k))$.
(2)~The oriented pretzel surface $P(t_1,\dots,t_k)$ is
bounded by a knot iff $k$ and $t_1,\dots,t_k$ are all odd.  In 
this case $g_s(\Bd P(t_1,\dots,t_k))$ is readily bounded below
using 5.3 (and the first remark).  In particular, it is easy to
recover the non-sliceness results of Yu and Kuga 
for certain algebraically slice pretzel knots \cite{12}
and connected sums thereof \cite{1}, obtained by methods 
of ``classical'' gauge theory 
(in contrast to the methods used to prove the local Thom Conjecture, 
an appeal to the truth of which underlies \cite{7} and thus the 
present paper).
\endremark

\enddocument


\Refs\widestnumber\no{10}
%
\ref\no 0
\by Michel Boileau and Lee Rudolph
\paper $\C^3$-algebraic Stein fillings via branched covers and plumbing
\paperinfo in preparation
\endref
%
\ref\no 1
\by Ken'ichi Kuga
\paper Making some non-slice knots
\jour J. College Arts Sci. Chiba Univ. B 
\vol 26 
\yr 1993
\pages 1--5
\endref
%
\ref\no 2
\by Takuji Nakamura 
\paper Four-genus and unknotting number of positive knots and links
\jour Osaka J. Math
\toappear
\endref
%
\ref\no 3
\by Lee Rudolph
\paper Braided surfaces and Seifert ribbons for closed braids
\jour Comment. Math. Helvetici
\vol 58
\yr 1983
\pages 1--37
\endref
%
\ref\no 4
\bysame
\paper Constructions of quasipositive knots and links, \rom{I}
\inbook Noeuds, Tresses, et Singularit\'es 
\yr 1983
\ed C. Weber
\publaddr Kundig, Geneva
\pages 233--246
\endref
%
\ref\no 5
\bysame
\paper A characterization of quasipositive Seifert surfaces 
\rom(Constructions
of quasipositive knots and links, \rom{III)}
\jour Topology
\vol 31
\yr 1992
\pages 231--237
\endref
%
\ref\no 6
\bysame
\paper Quasipositive annuli
\rom{(}Constructions of quasipositive knots and links, \rom{IV)} 
\jour J. Knot Theory Ramif. 
\vol 1
\yr 1993
\pages 451-466
\endref
%
\ref\no7
\bysame
\paper Quasipositivity as an obstruction to sliceness 
\jour Bull. A.M.S. 
\vol 29 
\yr 1993 
\pages 51--59
\endref
%
\ref\no 8
\bysame
\paper Quasipositive plumbing
\rom{(}Constructions of quasipositive knots and links, \rom{V)} 
\jour Proc. Amer. Math. Soc. 
\vol 126
\yr 1998
\pages 257--267
\endref
%
\ref\no9
\bysame
\paper Positive links are strongly quasipositive
\inbook Proceedings of the Kirbyfest 
\bookinfo Geometry and Topology Monographs 
\vol 2
\eds J. Hass and M. Scharlemann 
\yr 1999
\pages 555--562
\finalinfo math.GT/9804003
\endref
%
\ref\no10
\bysame
\paper Hopf plumbing, 
arborescent Seifert surfaces, 
baskets, espaliers, and homogeneous braids
\paperinfo preprint (1998)
\finalinfo math.GT/9810095
\endref
%
\ref\no11\by H. Seifert
\paper \"Uber das Geschlecht von Knoten
\jour Math. Annalen
\vol 110
\yr 1934
\pages 571--592
\endref
%
\ref\no 12\by Bao Zhen Yu
\paper A note on an invariant of Fintushel and Stern
\jour Topology and its Applications
\vol 38 \yr 1991
\pages 137--145
\endref
%
\endRefs
\bye